\documentclass[11pt]{amsart}
\usepackage{amsmath,amsfonts,latexsym,amssymb,amscd, graphicx,pictexwd}

\renewcommand{\P}{{\mathbb P}}
\newcommand{\E}{{\mathbb E}}
\newcommand{\R}{{\mathbb R}}
\newcommand{\LL}{{\mathbb L}}

\newcommand{\PP}{{\mathbf P}}

\newcommand{\0}{{\mathbf 0}}
\newcommand{\pp}{{\mathbf p}}
\newcommand{\qq}{{\mathbf q}}

\newcommand{\nn}{{\mathbf n}}
\newcommand{\cal}{\mathcal}

\newcommand{\calN}{{\mathcal N}}

\newcommand{\Z}{{\mathbb Z}}

\newcommand{\Var}{{\rm Var}}

\newtheorem{thm}{Theorem}
\newtheorem{coro}{Corollary}

\newtheorem{rem}{Remark}

\begin{document}
\title[Influence of the initial condition in equilibrium LPP models]{Influence of the initial condition in equilibrium last-passage percolation models}

\author{Eric Cator}
\address{Delft University of Technology\\
Mekelweg 4, 2628 CD Delft, The Netherlands}
\email{E.A.Cator@tudelft.nl}

\author{Leandro P. R. Pimentel}
\address{Universidade Federal do Rio de Janeiro\\
Caixa Postal 68530, 21941-909 Rio de Janeiro, RJ, Brasil}
\email{leandro@im.ufrj.br}

\author{Marcio Watanabe}
\address{Universidade de S\~ao Paulo\\
Caixa Postal 66281, 05311-970 S\~ao Paulo, SP, Brasil }
\email{msouza@ime.usp.br}


\begin{abstract}
In this paper we consider an equilibrium last-passage percolation model on an environment given by a compound two-dimensional Poisson process. We prove an $\LL^2$-formula relating the initial measure with the last-passage percolation time. This formula turns out to be a useful tool to analyze the fluctuations of the last-passage times along non-characteristic directions.
\end{abstract}

\maketitle

\section{Introduction and the main result}\label{sec:intro}
\subsection{The last-passage percolation model}
Let $\PP\subseteq\R^2$ be a two-dimensional Poisson random set of intensity one. On each point $\pp\in\PP$ we put a random positive weight $\omega_\pp$ and we assume that $\{\omega_\pp\,:\,\pp\in\PP\}$ is a collection of i.i.d. random variables, distributed according to a distribution function $F$, which are also independent of $\PP$. Throughout this paper we will make the following assumption on the distribution function $F$ of the weights:
\begin{equation}\label{eq:a2}
 \int_0^\infty e^{a x}\,dF(x) <+\infty\,,\mbox{ for some }a>0\,.
\end{equation}
This condition was used in \cite{CaPi} to prove the existence of invariant measures for the Hammersley's interacting fluid process we will introduce below. For each $\pp,\qq\in\R^2$, with $\pp<\qq$ (inequality in each coordinate, $\pp\neq \qq$), let $\Pi(\pp,\qq)$ denote the set of all increasing (or up-right) paths, consisting of points in $\PP$, from $\pp$ to $\qq$, where we exclude all points that share (at least) one coordinate with $\pp$. So we consider the points in the rectangle $]\pp,\qq]$, where we leave out the south and the west side of the rectangle. The  last-passage time between  $\pp\leq\qq$ is defined by
$$L(\pp,\qq):=\max_{\varpi\in\Pi(\pp,\qq)}\big\{\sum_{\pp'\in\varpi}\omega_{\pp'}\big\}\,.$$
When $F$ is the Dirac distribution concentrated on $1$ (each point has weight $1$ and we will denote this $F$ by $\delta_1$), then we refer to this model as the classical Hammersley model  \cite{AlDi,Ha}.

A crucial result is the following shape theorem (see Theorem 1.1 in \cite{CPshape}, p.164): set $\0=(0,0)$, $\nn=(n,n)$,
\begin{equation}\label{eq:constant}
\gamma=\gamma(F)=\sup_{n\geq 1}\frac{\E(L(\0,\nn))}{n}> 0\,\,\mbox { and }\,\,  f(x,t):=\gamma \sqrt{xt}\,.
\end{equation}
Then $\gamma(F)<\infty$ and for all $x,t>0$,
\begin{equation}\label{eq:shape}
\lim_{r\to\infty}\frac{L\left(\0,(rx,rt)\right)}{r}=\lim_{r\to\infty}\frac{\E L\left(\0,(rx,rt)\right)}{r} = f(x,t)\,.
\end{equation}

\subsection{The interacting fluid system formulation}
It is well known that the classical Hammersley model has a representation as an interacting particle system \cite{AlDi,Ha}. The general model has a similar description, although a better name might be an interacting fluid system. We start by restricting the compound Poisson process $\{\omega_\pp\,:\,\pp\in\PP\}$ to $\R\times\R_+$. To each measure $\nu$ on $\R$ we associate a non-decreasing process $\nu(\cdot)$ defined by
$$\nu(x)=\left\{\begin{array}{ll} \nu((0,x]) & \mbox{for } x\geq 0\\
-\nu((x,0]) & \mbox{for } x<0.\end{array}\right.$$
Let $\calN$ be the set of all positive, locally finite measures $\nu$ such that
$$\liminf_{y\to-\infty}\frac{\nu(y)}{y}>0\,.$$
We need this condition to define the evolution of the process, since otherwise all mass will be pulled to minus infinity. The Hammersley interacting fluid system $(M^{\nu}_t\,:\,t\geq 0)$ will be defined as a Markov process with values in $\calN$, as was done in \cite{CaPi}. Its evolution is defined as follows: if there is a Poisson point with weight $\omega$ at a point $(x_0,t)$, then $M^\nu_{t}(\{x_0\}) = M^\nu_{t-}(\{x_0\}) + \omega$, and for $x> x_0$,
$$M^\nu_t((x_0,x]) = (M^\nu_{t-}((x_0,x]) - \omega)_+\,.$$
Here, $M^\nu_{t-}$ is the ``mass distribution'' of the fluid at time $t$ if the Poisson point at $(x_0,t)$ would be removed. To the left of $x_0$ the measure does not change. In words, the Poisson point at $(x_0,t)$ moves a total mass $\omega$ to the left, to the point $x_0$, taking the mass from the first available fluid to the right of $x_0$. (See Figure \ref{fig:hamprocess} for a visualization, in case of atomic measures, of the process inside a space-time box.)
\begin{figure}[tb]
\begin{center}
\strut
$$
\beginpicture
\footnotesize
\setcoordinatesystem units <0.07\textwidth,0.05\textwidth>
  \setplotarea x from 0 to 10, y from 0 to 12
\multiput {\bf $\bigstar$} at
4 3
6 8
/


\multiput {$\bullet$} at
2 0
5 0
8 0
0 1.5
0 6
/
\put {$(0,0)$} at 0.0 -0.5
\put {$(x,0)$} at 10 -0.5
\put {$(0,t)$} at 0.0 10.5
\put {$(x,t)$} at 10 10.5
\put {$t/2$} at -0.5 5

\setlinear
\plot 0 10  10 10  10 0 0 0 0 10
/
\plot 10 8 8 8
/
\plot 8 0 8 8 6 8 6 10
/
\plot 8 0 8 3 5 3
/
\plot 5 0 5 3 4 3 4 6 0 6
/
\plot 8 6 3 6
/
\plot 2 0 2 6
/
\plot 2 1.5 0 1.5
/
\multiput {$\blacktriangle$} at
2 0.75
5 1.5
8 1.5
8 4.5
8 7
6 9
4 4.5
2 3.75
/
\multiput {$\blacktriangleleft$} at
1 1.5
1 6
3 6
6 6
6.5 3
7 8
9 8
4.5 3
/

{\large
\put {\bf 5} at 2 -0.3
\put {\bf 3} at 5 -0.3
\put {\bf 7} at 8 -0.3
\put {\bf 4} at -0.3 1.5
\put {\bf 6} at -0.3 6
\put {\bf 4} at 3.8 2.7
\put {\bf 7} at 5.8 7.7
}

\put {$4$} at 4.1 4.5
\put {$4$} at 4.45 3.3
\put {$3$} at 5.1 1.5
\put {$1$} at 6.45 3.3
\put {$7$} at 8.1 1.5
\put {$6$} at 8.1 4.5
\put {$5$} at 8.1 7
\put {$2$} at 8.95 8.3
\put {$7$} at 6.95 8.3
\put {$1$} at 5.95 6.3
\put {$5$} at 2.95 6.3
\put {$6$} at 0.95 6.3
\put {$1$} at 2.1 3.75
\put {$4$} at 0.95 1.8
\put {$5$} at 2.1 .75
\put {$7$} at 6.1 9

\setdots<2pt>

\plot 0 5 10 5
/

\endpicture
$$
\caption{In this picture,  restricted to $[0,x]$, the measure $\nu$ consists of three atoms of weight $5$, $3$ and $7$. The Poisson process, restricted to $[0,x]\times[0,t]$, has two points with weights $4$ and $7$. The measure $M^\nu_{t/2}$ consists of three atoms of weight $1$, $4$ and $6$, while at time $t$, it consists of one atom with weight $7$. A total weight of $4+6$ has left the box due to Poisson points to the left of the box, while a total weight of $2$ has entered.}\label{fig:hamprocess}
\end{center}
\end{figure}

In this paper we follow the Aldous and Diaconis \cite{AlDi} graphical representation in the last-passage model (compare to the result in the classical case, found in their paper): For each $\nu\in\calN$, $x\in\R$ and $t\geq 0$ let
\begin{equation}\label{eq:Lnu def}
L_{\nu}(x,t):=\sup_{z\leq x} \left\{ \nu(z) + L((z,0),(x,t))\right\}\,.
\end{equation}
The measure $M_t^\nu$ defined by
$$M_t^\nu((x,y]):=L_\nu(y,t)-L_\nu(x,t)\,\mbox{ for }x<y\,,$$
defines a Markov process on $\calN$ and it evolves according to the Hammersley interacting fluid system \cite{CaPi}.

We now make the following important observation for a random initial condition $\nu$, which basically follows from translation invariance.
\begin{thm}\label{thm:transinv}
Suppose $\nu\in\calN$ is a random initial measure on $\R$ independent of the Poisson process in $\R\times \R^+$, whose distribution is translation invariant. For any speed $V\in\R$ and any $x\in \R$, we have
\[ L_\nu(Vt,t) \stackrel{\cal D}{=} L_\nu(x,t) - \nu(x-Vt).\]
\end{thm}

The relevance of this result is most clear when we consider equilibrium measures of the Hammersley's interacting fluid process. Assume that we have a probability measure defined on $\calN$ and consider $\nu\in \calN$ as a realization of this probability measure. We say that $\nu$ is time invariant for the Hammersley interacting fluid process (in law) if
$$M^\nu_t\stackrel{\cal D}{=}M^\nu_0=\nu\,\,\,\mbox{ for all }\,\,\,t\geq 0\,.$$
In this case, we also say that the underlying probability measure on $\calN$ is an equilibrium measure. It is known that there is only one family of ergodic equilibrium measures for the Hammersley interacting fluid system \cite{CaPi}. Let us denote it by $\{\nu_\lambda\,:\,\lambda>0\}$, where
\begin{equation}\label{eq:intensity}
\lambda:=\E\nu_\lambda(1)\,.
\end{equation}
For simple notation, put $L_\lambda:=L_{\nu_\lambda}$. The main result of this paper is the following formula:
\begin{coro}\label{thm:formula}
Recall \eqref{eq:constant} and \eqref{eq:intensity}, and let
\begin{equation}\label{eq:direction}
V_\lambda:=\left(\frac{\gamma}{2\lambda}\right)^2\,\,\mbox{ and }\,\, \psi_\lambda:=\frac{\gamma^2}{2\lambda}\,.
\end{equation}
Here, $V_\lambda$ is the characteristic speed corresponding to $L_\lambda$ and $\psi_\lambda$ is the growth rate of $L_\lambda(V_\lambda t,t)$. Then
\begin{equation}\label{eq:formula}
\E\big(\left\{L_\lambda(x,t)-\left[\nu_\lambda(x-V_\lambda t)+\psi_\lambda t\right]\right\}^2\big)=\Var\big(L_\lambda\left(V_\lambda t,t\right)\big)\,.
\end{equation}
\end{coro}

\subsection{A central limit theorem for the classical model} To illustrate the importance of \eqref{eq:formula}, let us restrict ourselves to the classical Hammersley model. In this set-up, the equilibrium measures are one-dimensional Poisson processes of intensity $\lambda$, and $\gamma=\gamma(\delta_1)=2$. Thus,
$$V_\lambda:=\frac{1}{\lambda^2}\,\,\mbox{ and }\,\, \psi_\lambda:=\frac{2}{\lambda}\,.$$
Cator and Groeneboom \cite{CaGr2} proved that the variance of $L_\lambda$ grows sub-linearly along the characteristic speed $\lambda^{-2}$. Together with Corollary \ref{thm:formula}, this implies
\begin{coro}\label{coro:dependence}
Let $(z_t)_{t\geq0}$ be a deterministic path.
Then
\begin{equation}\label{eq:dependence}
\lim_{t\to\infty}\frac{\E\big(\left\{L_\lambda(z_t,t)-\left[\nu_\lambda(z_t - \lambda^{-2}t)+2\lambda^{-1}t\right]\right\}^2\big)}{t}=\lim_{t\to\infty}\frac{\Var\big(L_\lambda\left(\lambda^{-2}t,t\right)\big)}{t}=0\,.
\end{equation}
\end{coro}

\noindent{\bf Proof of Corollary \ref{coro:dependence}:} Formula \eqref{eq:formula}, applied to the classical model, gives us
$$\E\big(\left\{L_\lambda(z_t,t)-\left[\nu_\lambda(x-\lambda^{-2}t)+2\lambda^{-1}t\right]\right\}^2\big)=\Var\big(L_\lambda\left(\lambda^{-2}t,t\right)\big)\,.$$
On the other hand, \cite{CaGr2} shows that
$$\lim_{t\to\infty}\frac{\Var\big(L_\lambda\left(\lambda^{-2}t,t\right)\big)}{t}=0\,,$$
which proves \eqref{eq:dependence}.
\hfill$\Box$\\

\begin{coro}\label{coro:clt}
Let $(z_t)_{t\geq0}$ be a deterministic path such that
$$\lim_{t\to\infty}\frac{z_t}{t}=a\,.$$
Then
\begin{equation}\label{eq:var}
\lim_{t\to\infty}\frac{\Var \big(L_\lambda(z_t,t)\big)}{t}=\sigma^2:=|a\lambda-\frac{1}{\lambda}|\,.
\end{equation}
Furthermore, if $a\neq \lambda^{-2}$ then
\begin{equation}\label{eq:clt}
\lim_{t\to\infty}\P\big(L_\lambda(z_t,t)\leq \lambda z_t+\frac{t}{\lambda} +(\sigma\sqrt{t}) u\big)=\P(N\leq u)\,,
\end{equation}
where $N$ is a standard Gaussian random variable.
\end{coro}

\noindent{\bf Proof of Corollary \ref{coro:clt}:} Corollary \ref{coro:dependence} shows that
$$\lim_{t\to\infty}\frac{L_\lambda(z_t,t)-\left[\nu_\lambda(z_t-\lambda^{-2}t)+2\lambda^{-1}t\right]}{\sqrt{t}}=0\,,$$
in the $\LL^2$ sense. Since $\nu_\lambda$ is a one-dimensional Poisson process of intensity $\lambda$, this implies \eqref{eq:var} and \eqref{eq:clt}.
\hfill$\Box$\\

\begin{rem}\label{rem:cuberoot} Cator and Groeneboom \cite{CaGr2} proved that $\sqrt{\Var\big(L_\lambda\left(\lambda^{-2}t,t\right)\big)}$ is of order $t^{1/3}$, which gives us the same order for the $\LL^2$-distance between $L_\lambda(z_t,t)$ and $\left[\nu_\lambda(x-\lambda^{-2}t)+2\lambda^{-1}t\right]$.
\end{rem}
\begin{rem}\label{rem:BaikRains}
The central limit theorem for $L_\lambda$ (along any direction) was proved by Baik and Rains \cite{BaRa}. Their method was based on very particular combinatorial properties of the classical model that do not seem to hold for the general set-up. Our approach reveals the strong relationship with the initial configuration.
\end{rem}

\begin{rem}\label{rem:universality}
In the general set-up, Corollary \ref{thm:formula} implies: If the variance of $L_\lambda$ along the characteristic speed $V_\lambda$ is sub-linear, and the equilibrium measure has Gaussian fluctuations, then $L_\lambda$ will also have Gaussian fluctuations along non-characteristic directions.
\end{rem}

\begin{rem}\label{rem:variance formula}
For the classical Hammersley process an important formula for the variance of $L_\lambda(x,t)$ was derived in \cite{CaGr2}, Theorem 2.1:
\[ \Var(L_\lambda(x,t)) = -\lambda x + \frac{t}{\lambda} + 2\lambda \E (x-X_\lambda(t))_+,\]
where $X_\lambda(t)$ is the position at time $t$ of a second class particle starting at zero.
This formula was pivotal in deriving the cube-root behavior of $L_\lambda$ in \cite{CaGr2}, and later corresponding formulas were used to prove cube-root behavior for TASEP \cite{BaCaSep} and for ASEP \cite{BaSep}. However, this formula does not directly show the relationship with the initial configuration. Also, there seems to be no direct way to deduce \eqref{eq:formula} from this formula, even if we reformulate it, as was done in Equation (3.6) of \cite{CaGr2}, in terms of the exit-point of the longest path from $(0,0)$ to $(x,t)$, which is the right-most $z$ for which the supremum in \eqref{eq:Lnu def} is attained.
\end{rem}

\section{Proof of Theorem \ref{thm:transinv} and Corollary \ref{thm:formula}}
Recall that
\[ L_\nu(x,t) = \sup_{z\leq x} \left\{ \nu(z) + L((z,0),(x,t))\right\}.\]
Clearly, $L((z,0),(Vt,t))\stackrel{\cal D}{=} L((z+x-Vt,0),(x,t))$. By assumption, $\nu$ has a translation invariant distribution, independent of $L$. This implies that $\{z\mapsto \nu(z)\}\stackrel{\cal D}{=} \{z\mapsto\nu(z+x-Vt)-\nu(x-Vt)\}$, and
\begin{eqnarray*}
L_\nu(Vt,t) & \stackrel{\cal D}{=} & \sup_{z\leq Vt} \left\{ \nu(z+x-Vt) - \nu(x-Vt) + L((z+x-Vt,0),(x,t))\right\}\\
& = & \sup_{z\leq x} \left\{ \nu(z) + L((z,0),(x,t))\right\} - \nu(x-Vt)\\
& = & L_\nu(x,t) - \nu(x-Vt).
\end{eqnarray*}
This proves Theorem \ref{thm:transinv}.\hfill $\Box$\\

Corollary \ref{thm:formula} now follows from results in \cite{CaPi}: there it is shown that for any speed $V$, the stationarity of $L_\lambda$ leads to
\[ \E L_\lambda(Vt,t) = V\lambda t + \frac14 \gamma^2 t/\lambda.\]
This follows from the fact that the Hammersley fluid process has intensity $\lambda$ on the bottom side of the rectangle between $(0,0)$ and $(x,t)$, and intensity $\gamma^2/(4\lambda)$ on the left side (this refers to the expected mass of the fluid leaving the interval $[0,x]$ through $0$ per time unit). When we define the characteristic speed $V_\lambda=\gamma^2/(4\lambda^{2})$, then
\[ \E L_\lambda(V_\lambda t,t) = \psi_\lambda t.\]
This together with Theorem \ref{thm:transinv} immediately shows \eqref{eq:formula}.\hfill$\Box$\\

\section{The lattice last-passage percolation model}
In the lattice last-passage percolation model one considers i.i.d. weights  $\{\omega_\pp\,:\,\pp\in\Z^2\}$, distributed according to a distribution function $F$. For $F(x)=1-e^{-x}$ (exponential weights), we have a similar shape theorem as \eqref{eq:shape} with limit shape given by
$$f(x,t)=(\sqrt{x}+\sqrt{t})^2\,.$$
We know from \cite{BaCaSep} that the invariant measures are given by
$$\nu_\rho((x,y])\stackrel{\cal D}{=}\sum_{z=x+1}^y X_z\,,$$
where $\{X_z\,:\,z\in\Z\}$ is a collection of i.i.d. exponential random variables with parameter $\rho$.
The analog to formula \eqref{eq:formula} is
$$\E\big(\left\{L_\rho(x,t)-\left[\nu_\rho(x-\lfloor V_\rho t\rfloor)+\psi_\rho t\right]\right\}^2\big)=\Var \big(L_\rho\left(\lfloor V_\rho t\rfloor,t\right)\big)\,,$$
where
$$V_\rho:=\frac{\rho^2}{(1-\rho)^2}\,\,\mbox{ and }\,\,\psi_\rho:=\frac{1}{(1-\rho)^2}\,.$$
Together with the cube-root asymptotics \cite{BaCaSep}, this implies that
$$\lim_{t\to\infty}\frac{\E\big(\left\{L_\rho(z_t,t)-\left[\nu_\rho(z_t-V_\rho t)+\psi_\rho t\right]\right\}^2\big)}{t}=0\,.$$
Therefore, if
$$\lim_{t\to\infty}\frac{z_t}{t}=a\,$$
then
$$\lim_{t\to\infty}\frac{\Var\big( L_\rho(z_t,t)\big)}{t}=\sigma^2:=\frac{|a(1-\rho)^2-\rho^2|}{\rho^2(1-\rho)^2}\,,$$
and if $a\neq V_\rho$ then
$$\lim_{t\to\infty}\P\left(L_\rho(z_t,t)\leq \frac{z_t}{\rho}+ \frac{t}{1-\rho} +(\sigma\sqrt{t}) u\right)=\P(N\leq u)\,,$$
where $N$ is a standard Gaussian random variable.

\begin{rem}
Ferrari and Fontes \cite{FeFo} determined the dependence on the initial condition for the totally asymmetric exclusion process, which is isomorphic to the lattice last-passage percolation model with exponential weights. The method developed in this paper resembles the ideas in their paper. Bal\'azs \cite{Ba} used a different method to get a generalization of the Ferrari-Fontes result for certain types of deposition models. It is not clear to us whether our methods would work for these more general deposition models.
\end{rem}

\begin{rem}
In the general lattice model, the shape theorem \eqref{eq:shape} holds. However, not much is known about the limit shape $f$. If this function would not be strictly curved (we know it is convex, so this would mean that there are ``flat'' pieces), then the methods used in \cite{CPshape} to prove the existence and uniqueness of semi-infinite geodesics in a fixed direction do not apply, and we are not able to prove the existence of equilibrium measures.
\end{rem}

\end{document}